\documentclass{llncs}
\usepackage{amsfonts,amssymb,amsmath}
\usepackage{graphicx}
\usepackage{color,soul}
\usepackage{epstopdf}
\usepackage{caption}

\begin{document}

\title{The inexact residual iteration method for quadratic eigenvalue problem and the analysis of convergence}
\titlerunning{ Inexact residual iteration method }  
%
\author{Liu Yang \and Yuquan Sun\inst{*} \and Fanghui Gong
 }
\authorrunning{Liu Yang et al.} 
%
\tocauthor{Liu Yang, Yuquan Sun and Fanghui Gong }
\institute{LMIB \& School of Mathematics and Systems Science, BeiHang University, Beijing China, 100191.\\
\inst{*}{corresponding author, sunyq@buaa.edu.cn}
}
\maketitle              

\begin{abstract}
In this paper, we first establish the convergence criteria of the residual iteration method for solving quadratic eigenvalue problems. We analyze the impact of shift point and the subspace expansion on the convergence of this method. In the process of expanding subspace, this method needs to solve a linear system at every step. For large scale problems in which the equations cannot be solved directly,  we propose an inner and outer iteration version of the residual iteration method. The new method uses the iterative method to solve the equations and uses the approximate solution to expand the subspace. We analyze the  relationship between inner and outer iterations and provide a quantitative criterion for the inner iteration which can ensure the convergence of the outer iteration. Finally, our numerical experiments provide proof of our analysis
 and demonstrate the effectiveness of the inexact residual iteration method.
\end{abstract}

\section{Introduction}

The quadratic eigenvalue problem (QEP) is to find scalars $\lambda$ and nonzero vectors $x,y$ satisfying
\begin{equation}\label{QEPO}
Q(\lambda)x=(\lambda^{2}M+\lambda C+K)x=0, x\neq 0,
\end{equation}
\begin{equation}
y^{*}Q(\lambda)=y^{*}(\lambda^{2}M+\lambda C+K)=0,y\neq 0,
\end{equation}
where $M,C,K$ are $n\times n$ complex matrices, $x,y$  are the right and left eigenvectors, respectively, corresponding to the eigenvalue $\lambda$.
It has many important applications including least squares problem with constraints, fluid mechanics, circuit simulation, and structural mechanics. In\cite{Tisseur01}, Tisseur systematically summarized and reviewed the QEP.

There are two major classes of numerical methods to solve large QEPs. The first method is to linearize the QEP into an equivalent  generalized eigenvalue problem (GEP)  such as\cite{Tisseur00}

\begin{equation}\label{GEP}
A
\begin{bmatrix}
\lambda x\\  x
\end{bmatrix}
=
\lambda B
\begin{bmatrix} \lambda x\\  x
\end{bmatrix},
\end{equation}
\\
where
$$A=
\begin{bmatrix} -C&-K\\ I&0
\end{bmatrix},
B=
\begin{bmatrix} M&0\\  0&I
\end{bmatrix}.$$
When $M$ is reversible, we can transform it to an equivalent standard eigenvalue problem (SEP)  \\
\begin{equation}
B^{-1}A
\begin{bmatrix}
\lambda x\\
 x
\end{bmatrix}
=
\lambda
\begin{bmatrix}
\lambda x\\
 x
\end{bmatrix}.
\end{equation}

Then we can obtain the eigenpair $(\lambda, x)$ from the  eigenpair of $(A,B)$ or $B^{-1}A$. There are also other linearization  forms \cite{Tisseur01,SaadY92}.

The biggest advantage of this approach is that it can use all the theoretical and numerical results of standard and generalized eigenvalue problems. There are many well developed methods available, for example, rational Krylov method \cite{RuheA98}, displacement inverse Arnoldi method\cite{SaadY92,NatarajanR92}, Jacobi-Davidson method\cite{Sleijpen96,Van96}.
The disadvantages of this method are also obvious.
First, the size of the linearization is twice of the original QEP. This results in a significant increase in computational and storage requirements. Second, this method encounters stability problems \cite{Higham07}.

The second class of methods to solve large QEP are direct projection methods. These methods project the large QEP directly to a subspace $U$  to obtain a smaller QEP.  They can preserve the structure iteration of the original problem and has better numerical stability. Some direct projection methods, are the residual iteration method\cite{NeumaierA85,HuitfldtJ90,Jia97}, and the Jacobi-Davidson method\cite{Sleijpen96,Van96,Van96News,Van99}.

The main difficulty with these methods is the lack of theoretical basis.
QEPs are an important type of nonlinear eigenvalue problems that are less familiar and less routinely solved than the SEP and the GEP. A $n$ dimension QEP can have $2n$ eigenvalues (when $A$ is singular, the number of the eigenvalue is less than $2n$) and eigenvectors. SEP has Schur decomposition form, GEP has generalized Schur decomposition form, but QEP does not have such forms. For the direct projection methods, the properties of the projection subspace $U$ has not been thoroughly analyzed. We are not very familiar with the information contained in the subspace. Therefore, there are no perfect convergence analysis of these algorithms. Another disadvantage is that most methods require matrix inversion at each iteration. While the convergence speed is fast, the computational costs are larger.

In this paper, We first analyze the convergence of the residual iteration method. We show the property of the subspace which is constructed by the residual iteration process. We have established the relationship between the subspace and the desired eigenvectors. We analyze the impact of shift selection and subspace expansion on the capacity of the subspace containing desired eigenvectors.
In the process of expanding subspace, this method needs to solve a linear system at every step. For large scale problems in which the equations cannot be solved directly, then we propose an inner and outer iteration version of the residual iteration method. The new method uses iterative method to solve the equations and uses the approximate solution to expand the subspace. We establish the  relationship between inner and outer iteration and give a quantitative criterion for inner iteration which can ensure the convergence of outer iteration.

The remainder of this paper is arranged as follows. In section 2, we introduce the quadratic residual iteration method and analyze the convergence of this method. In section 3, we propose an inexact residual iteration method and give a quantitative convergence analysis of the method.
In section 4, several numerical experiments are presented to verify the results in this paper.

Throughout the paper, we denote by $||\cdot||$ the 2-norm of a vector or matrix, by $I$ the identity
matrix with the order clear from the context, by the superscript $‘*’$ the conjugate transpose of
a vector or matrix. We measure the distance between a nonzero vector $x$ and a subspace $\mathcal{V}$ by
$$\sin\angle(\mathcal{V},x)=\frac{||(I-P_{V})x||}{||x||}=\frac{||V^{*}_{\perp}x||}{||x||}$$
where $P_{V}$ is the orthogonal projector onto $\mathcal{V}$ and the columns of $V_{\perp}$ form an orthonormal
basis of the orthogonal complement of $\mathcal{V}$.

\section{Convergence analysis of the quadratic residual iteration method}

In this section, we first introduce the quadratic residual iteration method. Then we analyze the convergence property of the method. This provides a theoretical basis to further improve this type of methods.

For the most original residual iteration method, we  can find inspiration from Newton iterative method.
If we transform the quadratic eigenvalue problem into an equation
$$F(\lambda, x) = Q(\lambda)x=0,$$ and solve it using the Newton's Method, then we can obtain the following method.

{\bf Method 1:} One step quadratic residual iteration method
\begin{enumerate}
\item[1:]  for $k=0,1,2,\dots,$ do:
\item[2:]\hspace{0.3cm} $y^{(k)}=Q(\lambda_{k})^{-1}Q^{'}(\lambda_{k})x^{(k)}$.
\item[3:]\hspace{0.3cm}   $x^{(k+1)}=y^{(k)}/e^{*}y^{(k)}$.
\item[4:]\hspace{0.3cm}  $\lambda_{k+1}=\lambda_{k}-1/e^{*}y^{(k)}$.
\item[5:] end for
\end{enumerate}\label{newton}

where $Q^{'}(\lambda_{k})= 2 \lambda_{k} M + C$.
This is a very elementary method which can only find one eigenvalue and eigenvector. It is similar to the Rayleigh quotient iteration method for SEP. Depending on the characteristics of the Newton iteration method, it can have a faster convergence rate when it begins to converge. But the convergence result is affected by the initial value $(\lambda_{0},x^{(0)})$.
Under certain conditions, this method may have error convergence. In \cite{NeumaierA85}, some deficiencies are addressed resulting in an improved method. In process Method 1, we need one matrix inversion and one more matrix-vector product $Q'(\lambda_{k})x^{(k)}$.

When we use the subspace  projection method, we use the residual iteration process to expand the subspace\cite{MeerbergenK01} as follows:

{\bf Method 2:} Subspace quadratic residual iteration method (QRI)
\begin{enumerate}
\item[1:] set shift $\sigma$ and initial vector  $v_{1}\in R^{n}$, $||v_{1}||=1$ and convergence tolerance $tol$.
\item[2:] for $k=1,2,\dots,$ do:
\item[3:]\hspace{0.3cm}  set $V_{k}=[v_{1},v_{2},\ldots,v_{k}]$.
\item[4:]\hspace{0.3cm} project the large QEP \eqref{QEPO} onto subspace $V_{k}$ and solve the smaller QEP \\
   $$\omega^{2}M_{k}z+\omega C_{k}z+K_{k}z=0$$
   to get the eigenvalue $\omega_i$ and eigenvector $z_i$.
\item[5:]\hspace{0.3cm}  compute Ritz vector $\widetilde{x_{i}}=V_{i}z_{i}$ as approximate eigenvector.
\item[6:]\hspace{0.3cm}  compute the residual  $r_{i}=\omega^{2}_{i}M\widetilde{x_{i}}+\omega_{i}C\widetilde{x_{i}}+K\widetilde{x_{i}}$.
\item[7:]\hspace{0.3cm} select the first non-convergence  $r$.
\item[8:]\hspace{0.3cm}  compute $u=(\sigma^{2}M +\sigma C+K)^{-1}r$.
\item[9:]\hspace{0.3cm} get $v_{k+1}$ from $u$.
\item[10:] end for
\end{enumerate}

This method projects the large QEP to the subspace directly at step 4. It constructs the subspace by using the residual iteration method to generate new vector $v_{k+1}$ at each step. At step 7, if the first residual $\|r_1\|< tol$, we can use the second residual to expand the subspace. So we can compute more than one eigenpairs using this method. In Method 1, the new  vector is $u=Q(\omega_{i})^{-1}Q'(\omega_{i})r_{i}$. If we use it to expand the subspace in Method 2, it needs one more matrix-vector product. We can remove the term  $Q'(\omega_{i})$ and use a fixed value $\sigma$ in Method 2.
In order to get an orthogonal basis of the projection subspace, we should orthogonalize  $u$ with $v_{1},\dots,v_{k}$ at step 9.
At step 8, when the matrix size is small, the expanding vector $u$ can be computed directly. Usually, $u$ can be  obtained by solving the linear equations $u=(\sigma^{2}M +\sigma C+K)^{-1}r $.

Jia and Sun proposed a refined residual iteration method \cite{jia04}. They introduced the idea of refined projection to quadratic eigenvalue problem. For computing the approximate eigenvector of Ritz value $\omega$ at step 5, the refined vector $u$ satisfies
\begin{equation}
\|(\omega ^2 M + \omega C+ K) {u}\| =  \min_{
 \mbox{\scriptsize $u_i\in V_k$}}
\|(\omega ^2 M + \omega ^2 C+ K )u_i\|. \label{th2}
\end{equation}

This method also has some changes and modifications, for example,
we can use different shift $\sigma$ at each step. We can use vector $q =x_{k}-Q(\sigma)^{-1}r_k$ to expand subspace. So this method has a close relationship with Jacobi-Davidson method \cite{Sleijpen96,Van96,Van96News,Van99}. All methods have shown good convergence properties. However the convergence of these methods have not thoroughly analyzed.
According to the property of the residual iteration method for SEP, with this method, it is easy to determine the eigenvalues which are close to the target point $\sigma$. So, beneficial shift can accelerate the convergence of the method. For a given target point $\sigma$, $\lambda$ is the closest eigenvalue and $x$ is the corresponding eigenvector. Then $(\lambda,x)$ is the desired eigenpair of the subspace $\mathcal{V}_{k}$.
We can use the distance between $x$ and $\mathcal{V}_{k}$ or the angle $\angle(\mathcal{V}_{k},x)$ to measure the convergence of the method.

Let $\mathcal{V}_{k}$ be a subspace produced by Method 2. We use $\mathcal{V}_{k}$ to denote the projection subspace and its orthogonal base.
Let $P_{V_{k}}$ be an orthogonal projection operator onto subspace $\mathcal{V}_{k}$. When we get $u=Q^{-1}(\sigma)r_{k}$ with $||u||=1$ and $r_{k}$ is residual,  define $V_{k+1}= [V_k,v_{k+1}]$, the expanded subspace can be written as $\mathcal{V}_{k+1} = span\{V_{k+1}\}$, where $v_{k+1}=(I-P_{V_{k}})u$.

To improve subspace expansion, we have the following demonstration:

{\bf Theorem 1}
Let $\mathcal{V}_{k}$ and $\mathcal{V}_{k+1}=span \{V_k,v_{k+1}\}$ where the subspace is computed by Method 2,  $(\lambda,x)$ is the desired eigenpair. Suppose $\sin\angle(v_{k+1},x_{\perp})\neq0$, and $x_{\perp}=(I-P_{V_{k}})x$, then we have：
\begin{equation}
\sin\angle(\mathcal{V}_{k+1},x)=\sin\angle(\mathcal{V}_{k},x)\sin\angle(v_{k+1},x_{\perp}).
\end{equation}

\begin{proof}
Because
$$\sin^{2}\angle(\mathcal{V}_{k},x)-\sin^{2}\angle(\mathcal{V}_{k+1},x)=||(I-P_{V_{k}})x||^{2}-||(I-P_{V_{k+1}})x||^{2}=|v^{*}_{k+1}x|^{2},$$ and $||x_{\perp}||=\sin\angle(\mathcal{V}_{k},x)$.
Then we have
\begin{equation}
\begin{aligned}
\frac{\sin\angle(\mathcal{V}_{k+1},x)}{\sin\angle(\mathcal{V}_{k},x)}
&=\sqrt{1-(\frac{|v^{*}_{k+1}x|}
  {\sin\angle(\mathcal{V}_{k},x)})^{2}}\\
&=\sqrt{1-(\frac{|v^{*}_{k+1}x_{\perp}|}{\sin\angle(\mathcal{V}_{k},x)})^{2}}\\
&=\sqrt{1-(\frac{||x_{\perp}||\cos\angle(v_{k+1},x_{\perp})}{\sin\angle(\mathcal{V}_{k},x_{\perp})})^{2}}\\
&=\sqrt{1-\cos^{2}\angle(v_{k+1},x_{\perp})}\\
&=\sin\angle(v_{k+1},x_{\perp}).
\end{aligned}
\end{equation}
Finally, we get
\begin{displaymath}
\sin\angle(\mathcal{V}_{k+1},x)=\sin\angle(\mathcal{V}_{k},x)\sin\angle(v_{k+1},x_{\perp}).
\end{displaymath}
\end{proof}

Since $x=P_{V_k}x+(I-P_{V_k})x=x_k+x_{\perp}$,
so we know that the value $\sin\angle(v_{k+1},x_{\perp})$ represents the ability of the method to obtain information from the outside subspace. At the $i$th step, we set $x_{i,\perp}=(I-P_{V_{i}})x$ then we can get the following result
\begin{equation}\label{multi2}
 \sin\angle(\mathcal{V}_{m},x)=  \prod^{m}_{i=1}\sin\angle(v_{i},x_{i,\perp})，
 \end{equation}
and
\begin{displaymath}
\sin\angle(v_{k+1},x_{\perp})=\min_{\beta\in R}||V_{k+1}-\beta x_{\perp}||.
\end{displaymath}

This means that the subspace $\mathcal{V}_{m}$ can contain more information of the desired eigenvector through expanding.
The value $\sin\angle(v_{k+1},x_{\perp})$  determines the effect of each extension. So it should be as small as possible and we can use the upper bound to estimate the convergence rate. To give this bound, we first need the following lemma.

Assume that $A-\lambda B$ is a generalized linear form of $Q(\lambda)$, then we have  \cite{Tisseur01}
\begin{equation}
Q(\lambda)^{-1}
 =-\begin{bmatrix}
I & 0
\end{bmatrix}
(A-\lambda B)^{-1}
\begin{bmatrix}
I\\
0
\end{bmatrix}.
   \end{equation}

{\bf Lemma\cite{jia04}} Let the diagonal elements of $\Lambda= diag(\lambda_{1},\dots,\lambda_{2n})$ are the eigenvalues of $Q(\lambda)$. $X=[x_{1},\dots,x_{2n}]$, $Y=[y_{1},\dots,y_{2n}]$ are the corresponding right and left eigenvectors. Assume that $\lambda $ is not an eigenvalue of $(A,B)$ and $Q(\lambda)$, then we have
\begin{equation}
Q(\lambda)^{-1}=X(\lambda I-\Lambda)^{-1}Y^{*}=\sum^{2n}_{i=1}\frac{x_{i}y^{*}_{i}}{\lambda-\lambda_{i}}.
   \end{equation}

For the convergence rate of the Method 2, we have the following result.

{\bf Theorem 2}
Let $v_{k+1}$, $\lambda_{i}$, $x_{i}$, $y_{i}$,($i=1,2,\dots,2n$) be the same as the earlier definitions. For a given shift $\sigma$,  $\lambda_{1}$ is the closest eigenvalue and $x_{1}$ is the corresponding eigenvector, $r$ is the residual, and $$|\lambda_{1}-\sigma|\ll|\lambda_{2}-\sigma|.$$
Let $$|\frac{1}{\lambda_{2}-\sigma}|=\max{|\frac{1}{\lambda_{i}-\sigma}|},
\xi=\frac{\sum^{2n}_{i=2}|y^{*}_{i}r|}{|y^{*}_{1}r|}, i=2,\dots, 2n. $$
Then we have
\begin{equation}
 \sin\angle(v_{k+1},x_{1,\perp})\leq|\frac{\lambda_{1}-\sigma}{\lambda_{2}-\sigma}|\cdot\xi .
 \end{equation}\label{th2.1}

\begin{proof}
From
\begin{displaymath}
Q^{-1}(\sigma)=\frac{x_{1}y^{*}_{1}}{\lambda_{1}-\sigma}+\frac{x_{2}y^{*}_{2}}{\lambda_{2}-\sigma}+\sum^{2n}_{i=3}\frac{x_{i}y^{*}_{i}}{\lambda_{i}-\sigma},
\end{displaymath}

we know
\begin{displaymath}
u=Q^{-1}(\sigma)r=\frac{x_{1}}{\lambda_{1}-\sigma}\cdot y^{*}_{1}r+\frac{x_{2}}{\lambda_{2}-\sigma}\cdot y^{*}_{2}r+\sum^{2n}_{i=3}\frac{x_{i}}{\lambda_{i}-\sigma}\cdot y^{*}_{i}r,
\end{displaymath}
then
\begin{displaymath}
(I-P_{V_{k}})u=(I-P_{V_{k}})x_{1}\cdot \frac{y^{*}_{1}r}{\lambda_{1}-\sigma}+(I-P_{V_{k}})x_{2}\cdot \frac{y^{*}_{2}r}{\lambda_{2}-\sigma}+(I-P_{V_{k}})\sum^{2n}_{i=3}x_{i}\cdot \frac{y^{*}_{i}r}{\lambda_{i}-\sigma},
\end{displaymath}

so
\begin{displaymath}
\frac{\lambda_{1}-\sigma}{y^{*}_{1}r}\cdot (I-P_{V_{k}})u=(I-P_{V_{k}})x_{1}+(I-P_{V_{k}})\frac{\lambda_{1}-\sigma}{y^{*}_{1}r}\cdot \sum^{2n}_{i=2}x_{i} \frac{y^{*}_{i}r}{\lambda_{i}-\sigma}.
\end{displaymath}
Since $v_{k+1}=(I-P_{V_{k}})u$, $(I-P_{V_{k}})x_{1}=x_{1,\perp}$, then it holds that\\
\begin{equation}
\begin{aligned}
||\frac{\lambda_{1}-\sigma}{y^{*}_{1}r}\cdot v_{k+1}-x_{1,\perp}||
&=|\frac{\lambda_{1}-\sigma}{y^{*}_{1}r}|\cdot||\sum^{2n}_{i=2}x_{i} \frac{y^{*}_{i}r}{\lambda_{i}-\sigma}||\\
&\leq |\frac{\lambda_{1}-\sigma}{y^{*}_{1}r}|\cdot\sum^{2n}_{i=2}|\frac{1}{\lambda_{i}-\sigma}|\cdot|y^{*}_{i}r|\cdot||x_{i}||\\
&=|\frac{\lambda_{1}-\sigma}{y^{*}_{1}r}|\cdot\sum^{2n}_{i=2}|\frac{1}{\lambda_{i}-\sigma}|\cdot|y^{*}_{i}r|\\
&\leq|\frac{\lambda_{1}-\sigma}{y^{*}_{1}r}|\cdot|\frac{1}{\lambda_{2}-\sigma}|\cdot\sum^{2n}_{i=2}|y^{*}_{i}r|\\
&=|\frac{\lambda_{1}-\sigma}{\lambda_{2}-\sigma}|\cdot\frac{\sum^{2n}_{i=2}|y^{*}_{i}r|}{|y^{*}_{1}r|}\\
&=|\frac{\lambda_{1}-\sigma}{\lambda_{2}-\sigma}|\cdot\xi
\end{aligned}.
\end{equation}
Finally we get
\begin{equation}
\sin\angle(v_{k+1},x_{1,\perp})\leq|\frac{\lambda_{1}-\sigma}{\lambda_{2}-\sigma}|\cdot\xi.\label{ubound2}
\end{equation}
\end{proof}

We can use this result to show the convergence of both Method 1 and Method 2. When we use a
fixed shift $\sigma$ in Method 1, we get a constant matrix $Q(\sigma)=\sigma^2 M+\sigma C +K$ and Method 1 becomes to the power method of matrix $E = Q(\sigma)^{-1}$. Let $|\mu_1|> |\mu_2|\cdots\geq |\mu_n|$ are the eigenvalues of $E$. Then we can get the eiagenpair  $(\mu_1,q_1)$ by the power method.  If $\lambda_{1}$ is the nearest eigenvalue to
$\sigma$, $x_1$ is near to $q_1$ but there must be a constant gap between $x_1$ and $q_1$. We cannot find the exact eigenvector $x_1$ from constant matrix $E$. In Method 1, we use constantly changing shift $\sigma_k$ to get  changing  matrix $ Q(\sigma_k)$.

We have two ways to achieve the convergence. The first is to use better shift at each iteration. Method 1 is one such way.  At each iteration, we use new shift $\sigma^{(k)}$ to get new approximate eigenvector $x^{(k)}$. From  \eqref{ubound2}, we have
\begin{equation}
\sin\angle(x^{(k)},x_{1,\perp})\leq \left|\frac{\lambda_{1}-\sigma^{(k)}}{\lambda_{2}-\sigma^{(k)}}\right|\cdot\xi.\label{ubound3}
\end{equation}
When $\sigma^{(k)}$ is converging to $\lambda_{1}$,  $\left|\frac{\lambda_{1}-\sigma^{(k)}}{\lambda_{2}-\sigma^{(k)}}\right|\leq \left|\frac{\lambda_{1}-\sigma^{(k-1)}}{\lambda_{2}-\sigma^{(k-1)}}\right|$.
We can get better $x^{(k)}$ from $\sigma^{(k)}$, similarly,
better $x^{(k)}$ can give better $\sigma^{(k)}$.

The second way to achieve the convergence is subspace expanding. Now, we use fixed shift $\sigma$  in the iterative process and we obtain the better approximate eigenvector in the expanded subspace with factor $\left|\frac{\lambda_{1}-\sigma}{\lambda_{2}-\sigma}\right| \cdot\xi$.
Combining \eqref{multi2} and \eqref{ubound2}, we obtain
an estimation of total convergence rate for Method 2,
\begin{equation}
 \sin\angle(\mathcal{V}_{m},x) \leq (|\frac{\lambda_{1}-\sigma}{\lambda_{2}-\sigma}|)^m\cdot\xi^m .
 \end{equation}
Assume $\xi$ a moderate value at each iteration then
the convergence rate is decided by the value $|\frac{\lambda_{1}-\sigma}{\lambda_{2}-\sigma}|$.
The size of this value depends on the distribution of eigenvalues and selection of shift $\sigma$. So a good shift can accelerate the convergence of the method.

\section{Inexact quadratic residual iteration method}

The good convergence  of the quadratic residual iteration method mainly benefit from the special expanding vector $u$. When we compute the expanding  vector at step 8 of Method 2, we do not compute the inverse matrix directly. We can compute $u$ by solving  the corresponding linear equations. Usually, this is one of the most resource intensive  part of the method. For large scale problems, it is hard to compute the expanding vector, even by solving linear equations. For linear eigenvalue problems, one way to solve this problem is to compute an inexact solution of the linear equations. This kind of methods are called inexact method or outer inner iteration methods\cite{Jia97,jia04}.

Based on this principle, we propose the inexact iteration method for quadratic eigenvalue problems.

{\bf Method 3:} Inexact quadratic residual iteration method
\begin{enumerate}
\item[1:] set shift $\sigma$ and initial vector  $v_{1}\in R^{n}$, $||v_{1}||=1$ and convergence tolerance $tol$.
\item[2:] for $k=1,2,\dots,$ do:
\item[3:]\hspace{0.3cm}  set $V_{k}=[v_{1},v_{2},\ldots,v_{k}]$.
\item[4:]\hspace{0.3cm}   project the large QEP \eqref{QEPO} onto subspace $V_{k}$ and solve the smaller QEP \\
   $$\omega^{2}M_{k}z+\omega C_{k}z+K_{k}z=0$$
   to get the eigenvalue $\omega_i$ and eigenvector $z_i$.
\item[5:]\hspace{0.3cm}  compute Ritz vector $\widetilde{x_{i}}=V_{i}z_{i}$ as approximate eigenvector.
\item[6:]\hspace{0.3cm}  compute the residual  $r_{i}=\omega^{2}_{i}M\widetilde{x_{i}}+\omega_{i}C\widetilde{x_{i}}+K\widetilde{x_{i}}$.
\item[7:]\hspace{0.3cm} select the first non-convergence  $r$.
\item[8:]\hspace{0.3cm}solve the equations  $$(\sigma^{2}M +\sigma C+K)u=r.$$
\item[9:]\hspace{0.3cm} get $v_{k+1}$ from $u$.
\item[10:] end for
\end{enumerate}

The difference between Method 2 and Method 3 is at Step 8. Method 2 computes an exact solution, but Method 3 only needs to calculate an approximate solution. For large scale problems, it is difficult to solve the equations exactly. It is feasible to use a method to compute an approximate solution. Normally, we use an iterative method to compute the approximate solution of the equations. So, we call it inner iteration and call the iterative of Method 3 the outer iteration. The main difficulty of the method is to determine the accuracy of the approximate solution. We must to find the balance between outer and inner iteration. Let $u$ be the exact solution of and  $\widetilde{u}$ be the approximate solution. Then the relative error between them is
\begin{equation}\label{e3}
\varepsilon=\frac{||\widetilde{u}-u||}{||u||}.
\end{equation}
Then we can write
$$\widetilde{u}=u+\varepsilon||u||f,$$
where $f$ is the normalized error direction vector.\\
So we get:\\
\begin{equation}\label{e4}
(I-P_{V})\widetilde{u}=(I-P_{V})u+\varepsilon||u||f_{\perp},
\end{equation}
where
\begin{equation}\label{e5}
f_{\perp}=(I-P_{V})f.
\end{equation}

Define
\begin{equation}\label{e6}
\widetilde{v}=\frac{(I-P_{V})\widetilde{u}}{||(I-P_{V})\widetilde{u}||},
v=\frac{(I-P_{V})u}{||(I-P_{V})u||},
\end{equation}
and
\begin{equation}\label{e7}
\widetilde{\varepsilon}=\frac{||(I-P_{V})\widetilde{u}-(I-P_{V})u||}{||(I-P_{V})u||},
\end{equation}
where $\widetilde{v}$ and $v$ are the normalized subspace expansion vectors in the inexact and exact methods, respectively.
We can measure the difference between $\widetilde{v}$ and $v$ by
$\widetilde{\varepsilon}$ or $\sin\angle(\widetilde{v},v)$.
The relationship between $\widetilde{\varepsilon}$ and $\sin\angle(\widetilde{v},v)$ is \cite{jia14},
\begin{equation}\label{e8}
\sin\angle(\widetilde{v},v)=\widetilde{\varepsilon}\sin\angle(\widetilde{v},f_{\perp}).
\end{equation}
The relationship between $\widetilde{\varepsilon}$ and $\varepsilon$ is,
\begin{equation}
\varepsilon =\frac{||(I-P_{V})u||}{||u||\sin\angle(\mathcal{V},f)}\widetilde{\varepsilon}.
\end{equation}
Then we obtain:
\begin{equation}
\begin{aligned}
\widetilde{\varepsilon}=\frac{||u||\sin\angle(\mathcal{V},f)}{||(I-P_{V})u||} \varepsilon=\frac{\sin\angle(\mathcal{V},f)}{||(I-P_{V})\frac{u}{||u||}||} \varepsilon=\frac{\sin\angle(\mathcal{V},f)}{\sin\angle(\mathcal{V},u)}\varepsilon .
\end{aligned}
\end{equation}

For the standard eigenvalue problem, we have obtained the convergence properties of the exact method. We can prove the convergence of the inexact method by proving that $\widetilde{v}$ can mimic $v$.
According to the above relationships, we can show that
$\widetilde{v}$ can be a good imitation of $v$ under a moderately precise inner iteration.
But these convergence results are more or less based on prior knowledge of eigenvalues. For the quadratic eigenvalue problem, we need to analyze the requirements of the inner iteration directly from the convergence of the inexact method.

The convergence condition of Method 3 is $\sin\angle(\mathcal{V}_{m+1},x_{1})<\sin\angle(\mathcal{V}_{m},x_{1})$, $(\lambda_{1}, x_{1})$ is the desired eigenpair.
According to the conclusion $\sin\angle(\mathcal{V}_{m},x_{1})=\prod^{m}_{i=1}\sin\angle(v_{i},x_{i,\perp})$,
it is equivalent to
\begin{equation}
\sin\angle(\widetilde{v},x_{1,\perp})<1.
\end{equation}
According to their relationships of $\widetilde{v},\widetilde{u},
x_{1},x_{1,\perp}$, the convergence can be analyzed by the value $\sin\angle(\widetilde{u},x_{1})$ or $\tan\angle(\widetilde{u},x_{1})$.

For Method 3, we can take the approximate solution $\widetilde{u}$  as an exact solution of the following perturbed equation
\begin{equation}
 (\sigma^{2}M+\sigma C+K+\delta H)\widetilde{u} = r,
\end{equation}
here $\delta H$ is the perturbation matrix of $(\sigma^{2}M+\sigma C+K)$.

{\bf Lemma 1} If $||(\sigma^{2}M+\sigma C+K)^{-1}\delta H||<1$, the  approximate solution $\widetilde{u}$ and the exact solution $u$  have the following relationship
\begin{equation}\label{th3}
u-\widetilde{u} \approx (\sigma^{2}M+\sigma C+K)^{-1}\delta H u.
\end{equation}

\begin{proof}
For a matrix $X$ and the corresponding unit matrix $I$, if $||X||<1$, then $I-X$ is invertible\cite{J97} and
\begin{equation}\label{th1}
(I-X)^{-1}=\sum^{\infty}_{i=0}X^{i}.
\end{equation}

Now, we can get
\begin{equation}
\begin{aligned}
\widetilde{u}
&=(\sigma^{2}M+\sigma C+K+\delta H)^{-1} r\\
&=[I+(\sigma^{2}M+\sigma C+K)^{-1}\delta H]^{-1} (\sigma^{2}M+\sigma C+K)^{-1} r.
\end{aligned}
\end{equation}
We use formula\eqref{th1} to $[I+(\sigma^{2}M+\sigma C+K)^{-1}\delta H]^{-1}$ and ignore higher order terms to obtain
\begin{displaymath}
\begin{aligned}
\widetilde{u}
&\approx[I-(\sigma^{2}M+\sigma C+K)^{-1}\delta H] (\sigma^{2}M+\sigma C+K)^{-1} r\\
&=[I-(\sigma^{2}M+\sigma C+K)^{-1}\delta H] u\\
&=u-(\sigma^{2}M+\sigma C+K)^{-1}\delta H u.
\end{aligned}
\end{displaymath}
Then we obtain the relationship about $\widetilde{u}$ and $u$.
\end{proof}

For the convergence condition of Method 3, we have the following result.

{\bf Theorem 4} Suppose $(\lambda_{1}, x_{1})$ is the desired eigenpair, and $\tan\angle(u,x_{1})<\tan\angle(u-\widetilde{u},x_{1})$, then the angle between the inexact solution and the desired eigenvector satisfies
$$\tan\angle(u,x_{1})<\tan\angle(\widetilde{u},x_{1})<\tan\angle(u-\widetilde{u},x_{1}).$$
\begin{proof}
From
\begin{displaymath}
Q^{-1}(\sigma)=\frac{y^{*}_{1}}{\lambda_{1}-\sigma}x_{1}+\frac{y^{*}_{2}}{\lambda_{2}-\sigma}x_{2}+\sum^{2n}_{i=3}\frac{y^{*}_{i}}{\lambda_{i}-\sigma}x_{i},
\end{displaymath}
we can get
\begin{equation}\label{8}
u=Q^{-1}(\sigma)r=\frac{y^{*}_{1}r}{\lambda_{1}-\sigma}\cdot x_{1} +\frac{y^{*}_{2}r}{\lambda_{2}-\sigma}\cdot x_{2}+\sum^{n}_{i=3}\frac{y^{*}_{i}r}{\lambda_{i}-\sigma}\cdot x_{i}.
\end{equation}

When $|\lambda_{1}-\sigma|\ll|\lambda_{2}-\sigma|, |\frac{1}{\lambda_{2}-\sigma}|=\max{|\frac{1}{\lambda_{i}-\sigma}|}(i=2,\dots,2n)$, we can ignore the small value $\sum^{2n}_{i=3}\frac{1}{\lambda_{i}-\sigma}$. Then the formula \eqref{8} can be written as
\begin{equation}
u=Q^{-1}(\sigma)r=\frac{y^{*}_{1}r}{\lambda_{1}-\sigma}\cdot x_{1} +\frac{y^{*}_{2}r}{\lambda_{2}-\sigma}\cdot x_{2}.
\end{equation}
We can rewrite this formula as
\begin{equation}\label{q2}
u=\alpha_{1} x_{1} +\beta_{1} x_{1,\perp}.
\end{equation}
Similarly, the formula \eqref{th3} also can be written as
\begin{equation}\label{q3}
u-\widetilde{u}= \alpha_{2} x_{1} +\beta_{2} x_{1,\perp}.
\end{equation}
Therefore, combining the last relation with \eqref{q2} establishes
\begin{equation}\label{q4}
\widetilde{u}= (\alpha_{1}+\alpha_{2}) x_{1} +(\beta_{1}+\beta_{2}) x_{1,\perp}.
\end{equation}

Then we can get
\begin{displaymath}
\tan\angle(u,x_{1})=\frac{\beta_{1}}{\alpha_{1}},\tan\angle(\widetilde{u},x_{1})=\frac{\beta_{1}+\beta_{2}}{\alpha_{1}+\alpha_{2}},
\tan\angle(u-\widetilde{u},x_{1})=\frac{\beta_{2}}{\alpha_{2}}.
\end{displaymath}

From the following relationship
\begin{equation}\label{}
\frac{\beta_{1}}{\alpha_{1}}\leq \frac{\beta_{1}+\beta_{2}}{\alpha_{1}+\alpha_{2}} \leq \frac{\beta_{2}}{\alpha_{2  }  },
\end{equation}
we can finish the proof.
\end{proof}

From Theorem 4, we know that the inexact method can have fast convergence when the difference between $\frac{\beta_{1}}{\alpha_{1}}$ and $\frac{\beta_{2}}{\alpha_{2}}$ is very small. Because $\frac{\beta_{1}}{\alpha_{1}} =O(\frac{\lambda_{2}-\sigma}{\lambda_{1}-\sigma})$  and
 $\frac{\beta_{2}}{\alpha_{2}}=O(\frac{\lambda_{2}-\sigma}{\lambda_{1}-\sigma})$, these two values are  always close, regardless of whether they are large or small.

When $\frac{\beta_{1}}{\alpha_{1}}$ is a very small value, it means that both angles between $u$ and $x$, $v$ and $x_{\perp}$ are small. So the exact residual iteration method can have fast convergence rate. For the inexact method, if the inner iteration is not very precise, there is a large error between  $\widetilde{u}$ and $u$.  At first glance, the convergence rate of the inexact method may not be very fast. But small $\frac{\beta_{1}}{\alpha_{1}}$ means small $\frac{\beta_{2}}{\alpha_{2}}$. So the angle between $\widetilde{u}$ and $x$ is also small. That is to say, the inexact method can have fast convergence rate in spite of the moderate accuracy of inner iteration.

When the value $\frac{\beta_{1}}{\alpha_{1}}$ is not too small, this means that the inexact method cannot get fast convergence speed using high accuracy in inner iteration. But the convergence of the inexact method can be guaranteed by $\sin\angle(\widetilde{v},x_{\perp})<q$, where $q$ is a constant value less than one. This condition is relatively easy to be satisfied and is independent of the inner iteration. The analysis shows that the convergence of the inexact method is mainly determined by the shift selection and subspace expansion.
However, there are limitations to improving the precision of the inner iteration.

\section{Numerical experiments}

Several numerical experiments are presented in this section to demonstrate the effectiveness of Method 3 and the analysis results.
For all examples, DIM stands for the dimension of matrix, TOL denotes the convergent precision of outer iteration, tol is the  precision of inner iteration. We use GMRES to solve the inner iteration. We show the CPU time of each part of the method and unit is in seconds. TOTAL expresses the total CPU time, CGMRES stands for the steps of gmres, TGMRES stands for the CPU time of GMRES, unit is in seconds. ITER is the number of iteration.
We select three different tol in our experiments, which respective denoted by  "inexact1", "inexact2", "inexact3".
In following figures, the horizontal axis is the dimension of subspace
and the vertical axis is the relative largest of the six residuals norm at each subspace expansion.

{\bf Example 1} For a fixed shift $\sigma$, let $Q(\sigma)=\sigma^{2}M+\sigma C+K$ and $E=Q(\sigma)^{-1}$. We show that the eigenvector of $E$ is different from the eigenvector of the corresponding quadratic eigenvalue problem.  So we cannot use  the power method of $E$ to compute the eigenvector of $Q(\lambda)$.

We use example in\cite{Tisseur01}. The matrices are
$M =
\begin{bmatrix}
0 & 6 & 0\\
0 & 6 & 0\\
0 & 0 & 1\\
\end{bmatrix}$,
$C =
\begin{bmatrix}
1 & -6 & 0\\
2 & -7 & 0\\
0 & 0 & 0\\
\end{bmatrix}$,
and $ K$ is a unite matrix.
The six eigenpairs $(\lambda_{k}, x_{k}), k = 1:6$ are

\begin{center}
\begin{tabular}{cccccccccccccc}\hline
k              &       1       &    2           &  3  &   4   &   5  &  6       &     \\ \hline
$\lambda_{k}$  & $\frac{1}{3}$ & $\frac{1}{2}$  &  1  &  $i$  & $-i$ & $\infty$ &  \\
$x_{k}$        &
$\begin{bmatrix}
1 \\
1\\
0\\
\end{bmatrix}$
&$\begin{bmatrix}
1 \\
1\\
0\\
\end{bmatrix}$
&$\begin{bmatrix}
0 \\
1\\
0\\
\end{bmatrix}$
&$\begin{bmatrix}
0 \\
0\\
1\\
\end{bmatrix}$
&$\begin{bmatrix}
0 \\
0\\
1\\
\end{bmatrix}$
&$\begin{bmatrix}
1 \\
0\\
0\\
\end{bmatrix}$&\\
\hline
\end{tabular}
\end{center}

 We choose $\sigma=0.9$, the eigenvalue which is closest to $\sigma $ is $\lambda=1$, the corresponding eigenvector is $(0,1,0)^{H}$. The eigenvalues of $E=Q(\sigma)^{-1}$ are
the eigenvalues of E are 0.735294117647062, 9.999999999999933, 0.552486187845304, and the eigenvector of the largest eigenvalue is
$\begin{bmatrix}
-0.287347885566346 \\
-0.957826285221151 \\
 0                \\
\end{bmatrix}$.

We can see that the corresponding eigenvector of  $E$ is
different from the wanted eigenvector.
Through the above example, we show that there is a difference between the eigenvectors of shifted standard eigenvalue problem and the quadratic eigenvalue problem.

{\bf Example 2} This example is arising from a finite element model of a linear spring in parallel with Maxwell elements (a Maxwell element is a spring in series with a dashpot) \cite{BetckeT11}. The quadratic matrix polynomial is $Q(\lambda) = \lambda^{2}M+\lambda C+K$, where the mass matrix $M$ is rank deficient and symmetric, the damping matrix $C$ is rank deficient and block diagonal, and the stiffness matrix $K$ is symmetric and has arrowhead structure.
The matrices have the form:
\begin{displaymath}
M=diag(\rho\widetilde{M}_{11},0),C=diag(0,\eta_{1}\widetilde{K}_{11},\dots,\eta_{m}\widetilde{K}_{m+1,m+1}),
\end{displaymath}
\begin{displaymath}
K=
\begin{bmatrix}
\alpha_{\rho}\widetilde{K}_{11} & -\xi_{1}\widetilde{K}_{12} & \cdots & -\xi_{m}\widetilde{K}_{1,m+1}\\
-\xi_{1}\widetilde{K}_{12}      &e_{1}\widetilde{K}_{22}     & 0 & 0 \\
\vdots  & 0 & \ddots & 0 \\
-\xi_{m}\widetilde{K}_{1,m+1} & 0 & 0 & e_{m}\widetilde{K}_{m+1,m+1}\\
\end{bmatrix},
\end{displaymath}
where $\widetilde{M}_{ij}$ and $\widetilde{K}_{ij}$ are element mass and stiffness matrices, $\xi_{i}$ and $e_{i}$ measure the spring stiffness, and $\rho$ is the material density. In this example, these matrices are randomly generated by the method of reference\cite{BetckeT11}. The size of the matrices is 20000.
We use the new method to compute 6 eigenvalues which are closest to the shift $\sigma=-1.8$. We select three precision $10^{-3}, 10^{-4},10^{-5}$ for the inner  iteration to show the effect of the precision of inner iteration on the outer iteration. The precision of outer iteration is $10^{-12}$. Table 1 and Figure 1 show the results.

\begin{figure}
\begin{center}
\resizebox*{13cm}{!}{\includegraphics{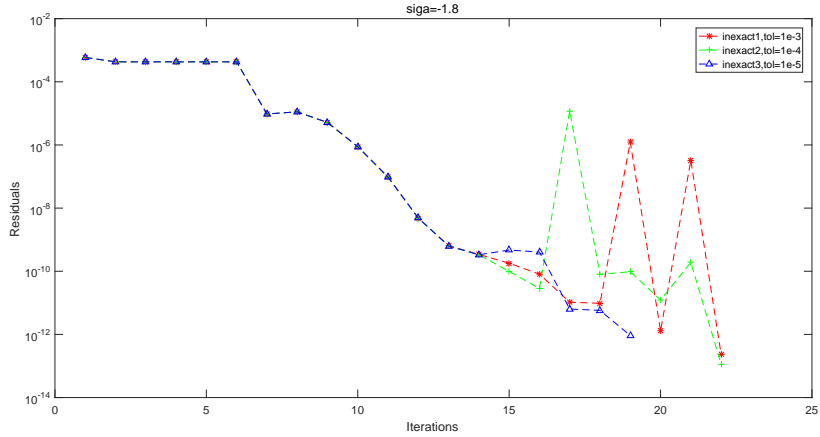}}
\captionsetup{labelsep=quad}
\caption{Example 2, $\sigma=-1.8$}
\label{pic1}
\end{center}
\end{figure}

\begin{table}[!h]
\begin{center}
\captionsetup{labelsep=quad}
\caption{Example 2, DIM=20000, TOL$={10}^{-12}$}
\label{table1}
\begin{tabular}{ccccccccccc}\hline
Standard & tol  & TOTAL & CGMRES &  TGMRES    &  ITER  &    \\ \hline
inexact1 & 1e-3 & 6.66 & 69 &1.41   &  22  &   \\
inexact2 & 1e-4 & 7.04 & 78 &1.48 &  22  &   \\
inexact3 & 1e-5 & 7.55 & 84 &1.54  &  19  &   \\
 \hline
\end{tabular}
\end{center}
\end{table}
It can be seen from Table 1 that the method with lowest precision of inner iteration uses more outer iterations than the method with highest precision of inner iteration, but the total time is less. Inner iteration using high precision, may reduce the number of external iteration but it need more number of inner iteration. So the inner iteration of low precision requires less time at each
outer iteration and has advantage on total CPU time.

Figure 1 shows the largest residual norm of different inner precision at each outer iteration.
They have the same convergence history in the first 14 outer iterations. The convergence of the higher  accuracy is more smooth than that of the lower.
But the method with the lowest precision can get the desired result using only three outer iterations.
In general, whether the CPU time or iteration number, the difference is not very obvious between methods with difference inner precision.

{\bf Example 3} This example was tested in \cite{BetckeT11}. It come from the finite element discretization 2D version of time-harmonic wave equation. On the unit square $[0,1]\times[0,1]$ with mesh size $h$,
the $n \times n$ coefficient matrices of $Q(\lambda)$ with $n=\frac{1}{h}(\frac{1}{h}-1)$ are given by:
\begin{displaymath}
M=-4\pi^{2}h^{2}I_{m-1}\otimes(I_{m}-\frac{1}{2}e_{m}e^{T}_{m}),C=2\pi i\frac{h}{\zeta}I_{m-1}\otimes (e_{m}e^{T}_{m}),
\end{displaymath}
\begin{displaymath}
K=I_{m-1}\otimes D_{m}+T_{m-1}\otimes(-I_{m}+\frac{1}{2}e_{m}e^{T}_{m}),
\end{displaymath}
where $\otimes$ denotes the Kronecker product, $m=\frac{1}{h}$, $\zeta$ is the (possibly complex) impedance, and
$D_{m}=\text{tridiag}(-1,4,-1)-2e_{m}e^{T}_{m}$, $T_{m-1}=\text{tridiag}(1,0,1)$.  We set $n=160000$,  and use the new method to compute 6  eigenvalues which are closest to the shift $\sigma=-0.5+4i$. We select three precision $10^{-3}, 10^{-4},10^{-5}$ for the inner iteration. The results are shown in Table 2 and Figure 2.

\begin{figure}[!h]
\begin{center}
\resizebox*{13cm}{!}{\includegraphics{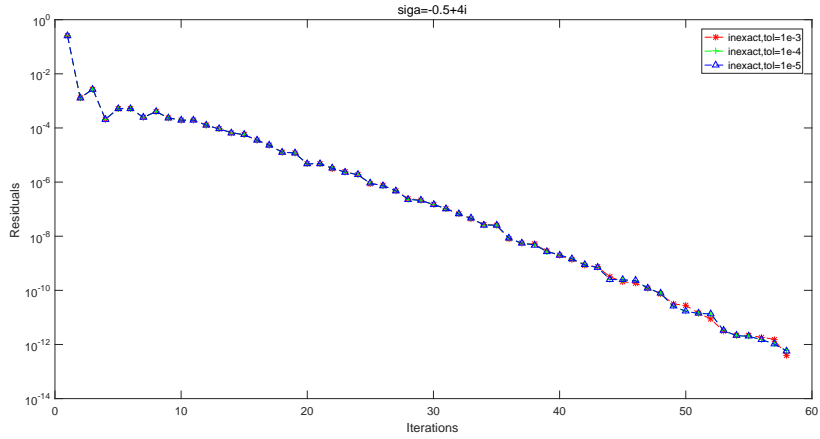}}
\captionsetup{labelsep=quad}
\caption{Example 3, $\sigma=-0.5+4i$}
\label{pic2}
\end{center}
\end{figure}

\begin{table}[!h]
\begin{center}
\captionsetup{labelsep=quad}
\caption{Example 3, DIM=160000, TOL=${10}^{-12}$}
\label{table 2}
\begin{tabular}{cccccccccccc}\hline
Standard & tol & TOTAL    & CGMRES &TGMRES   &  ITER  &    \\ \hline
inexact1  & 1e-3 & 663.27 & 471 &291.95 &  63  &   \\
inexact2  & 1e-4 & 972.78 & 672 &481.89   &  63  &   \\
inexact3  & 1e-5 & 1.1682e+03 &890 &599.73   &  63  &   \\
 \hline
\end{tabular}
\end{center}
\end{table}


For large scale problem, it needs more large subspace to find the desired eigenvalues. All of the three methods with different inner precision  need  more outer iteration. But, we can find out from Figure 2 that they have the same convergence history. There are no advantages if we use higher precision for inner iteration.

We can see from the table 2 that for the inner iteration with the lower precision, the less time it will use at each outer iteration. Then the total time saving is very considerable. So the method with lowest inner iteration is more superiority
than that using highest inner precision. This is consistent with our theoretical analysis.
In many practical applications, it is impossible,
even if we want to find a very accurate solution fot the inner iteration. So both theory and experiments show that this new method is a feasible and efficient method.

\newpage

\baselineskip 18pt
\renewcommand{\baselinestretch}{1.2}

\end{document}